\documentclass{elsart3-1}

\usepackage{amsmath,  amsopn, amsfonts}

\usepackage{amssymb}

\usepackage[english,francais]{babel}

\newtheorem{theorem}{Theorem}[section]
\newtheorem{lemma}[theorem]{Lemma}
\newtheorem{e-proposition}[theorem]{Proposition}
\newtheorem{corollary}[theorem]{Corollary}
\newtheorem{e-definition}[theorem]{Definition\rm}

\newtheorem{theoreme}{Th\'eor\`eme}[section]

\newtheorem{proposition}[theoreme]{Proposition}
\newtheorem{corollaire}[theoreme]{Corollaire}

\renewcommand{\theequation}{\arabic{equation}}
\def\og{\leavevmode\raise.3ex\hbox{$\scriptscriptstyle\langle\!\langle$~}}
\def\fg{\leavevmode\raise.3ex\hbox{~$\!\scriptscriptstyle\,\rangle\!\rangle$}}

\def\11{{\rm 1~\hspace{-1.4ex}l} }
\def\R{\mathbb R}

\def\cal{\mathcal}

\begin{document}

\begin{frontmatter}



\selectlanguage{english}
\title{Multilinear eigenfunction estimates for the Laplace spectral projectors on compact manifolds}

\vspace{-2.6cm}

\selectlanguage{francais}
\title{Estim{\'e}es multilin{\'e}aires pour les projecteurs spectraux du laplacien sur les vari{\'e}t{\'e}s compactes}

\author{N. Burq}
\ead{nicolas.burq@math.u-psud.fr}
\ead[url]{http://www.math.u-psud.fr/~burq}
\address{Universit{\'e} Paris Sud,
Math{\'e}matiques,
B{\^a}t 425, 91405 Orsay Cedex}
\author{P. G{\'e}rard}
 \ead{Patrick.gerard@math.u-psud.fr}
\address{Universit{\'e} Paris Sud,
Math{\'e}matiques,
B{\^a}t 425, 91405 Orsay Cedex}
\author{N. Tzvetkov}
\ead{nikolay.tzvetkov@math.u-psud.fr}
\ead[url]{http://www.math.u-psud.fr/~tzvetkov,}
\address{Universit{\'e} Paris Sud,
Math{\'e}matiques,
B{\^a}t 425, 91405 Orsay Cedex}
\selectlanguage{english}

\begin{abstract}
The purpose of this note is to extend to any space dimension the bilinear estimate for eigenfunctions of the Laplace operator on a compact manifold (without boundary) obtained in~\cite{BGT7} in dimension $2$. We also give some related trilinear estimates.
\vskip 0.5\baselineskip
\selectlanguage{francais}
\noindent{\bf R\'esum\'e}
\vskip 0.5\baselineskip
\noindent
L'objet de cette note est de g\'en\'eraliser {\`a} toute dimension d'espace les estimations bilin\'eaires de projecteurs spectraux de l'op{\'e}rateur de Laplace sur une vari{\'e}t{\'e} compacte (sans bord), d\'emontr\'ees dans~\cite{BGT7} en dimension 
~2. On {\'e}nonce aussi des estimations trilin\'eaires.
\end{abstract}
\end{frontmatter}

\selectlanguage{francais}
\section*{Version fran\c{c}aise abr\'eg\'ee}
Soit $(M,g)$ une vari{\'e}t{\'e} riemanienne compacte, $C^ \infty$ (sans bord) et $\Delta$ le laplacien sur les fonctions de $M$. Nous avons obtenu pr{\'e}c{\'e}demment (\cite{BGT7}) des estim{\'e}es bilin{\'e}aires sur les projecteurs spectraux du laplacien dans le cas o{\`u} la dimension de $M$ est $2$. Le but de cette note est de g{\'e}n{\'e}raliser ces estim{\'e}es en toute dimension d'espace:
\vskip 0.2 cm
\begin{theoreme}\label{th2}
Soit $\chi\in {\cal S}(\R)$. Pour $\lambda \in {\mathbb R}$ on note 
$\chi_\lambda = \chi(\sqrt{-\Delta}-\lambda)$ le projecteur spectral autour de $\lambda$.
Il existe $C$ tel que pour tous $\lambda, \mu \geq 1 $, $f,g \in L^{2}(M)$,
\begin{equation*}\label{eqbilin1}
\|\chi_\lambda f \, \chi_\mu g
\|_{L^{2}(M)}
\leq C \Lambda(d, \min (\lambda, \mu)) \|f\|_{L^{2}(M)}\|g\|_{L^{2}(M)}
,
\end{equation*}
{ avec }
\begin{equation*}
\Lambda(d,\nu) = \begin{cases}
  \nu^{\frac{1}{4}} &\text{ if $d=2$}\\
  \nu^{\frac{d-2}{2}}\log^{1/2} (\nu) &\text{ if $d=3$}\\
  \nu^{\frac{d-2}{2}} &\text{ if $d\geq 4$} 
\end{cases}
\end{equation*}
De plus on a aussi pour tous $\lambda,\mu,\nu \geq 1$, $f,g, h \in L^{2}(M)$ l'estimation trilin{\'e}aire suivante:
\begin{equation*}\label{eqtrilin1}
\|\chi_{\lambda}f\, \chi_{\mu} g \, \chi_{\nu}h
\|_{L^{2}(M)}
\leq
C\left(\frac{\lambda \mu \nu} { \max ( \lambda, \mu, \nu)}\right)^{\frac{2d-3}4}\|f\|_{L^{2}(M)}\|g\|_{L^{2}(M)}\|h\|_{L^{2}(M)}.
\end{equation*}
\end{theoreme}
Si on applique ce r{\'e}sultat en choisissant pour $f$ et $g$ deux harmoniques sph{\'e}riques sur la sph{\`e}re $\mathbb{S}^d$ on obtient
\vskip 0.2 cm
\begin{corollaire}
\label{coeigenf}
Il existe $C>0$ tel que si $H_p$ et $H_q$ sont deux harmoniques sph{\'e}riques de degr{\'e}s respectifs $p$ et $q$ plus grands que $1$, on a 
\begin{equation}\label{eqbilin1f}
\|H_p H_q
\|_{L^{2}(\mathbb{S}^d)}
\leq C \Lambda(d, \min(p,q)) \|H_p\|_{L^{2}(\mathbb{S}^d)}\|H_q\|_{L^{2}(\mathbb{S}^d)}
\end{equation}
\end{corollaire}
La version {\em lin{\'e}aire} de notre th{\'e}or{\`e}me ($\lambda= \mu=\nu$), sans la perte logarithmique pour $d=3$, est due {\`a} Sogge~\cite{SO86,SO88,S}. En dimension $2$, notre preuve dans~\cite{BGT7} {\'e}tait inspir{\'e}e par un travail de H{\"o}rmander \cite{Ho73} sur les op{\'e}rateurs satisfaisant {\`a} la condition de Carleson Sj{\"o}lin. Ici notre preuve est diff{\'e}rente (m{\^e}me pour $d=2$) et repose sur une bilin{\'e}arisation des arguments de~\cite{SO86,SO88,S} dans l'esprit de Stein~\cite{St93}. Plus pr{\'e}cis{\'e}ment, apr{\`e}s diverses r{\'e}ductions, on se ram{\`e}ne {\`a} {\'e}tablir deux estim{\'e}es (micro)-locales {\em lin{\'e}aires} sur l'op{\'e}rateur. La premi{\`e}re pr{\'e}cise la continuit{\'e} $L^2$ de l'op{\'e}rateur, tandis que la seconde est fond{\'e}e sur l'effet dispersif. Le fait que dans la relation~\eqref{eqbilin1f} la plus haute fr{\'e}quence disparaisse compl{\`e}tement du terme de droite {e}st crucial en vue des applications {\`a} l'{\'e}tude de l'{\'e}quation de Schr{\"o}dinger non lin{\'e}aire ({\em cf.}~\cite{BGT7}). Un calcul {\'e}l{\'e}mentaire sur les fonctions propres $e_n(x_1, x_2,x') = (x_1+ix_2)^n$, $(x_1, x_2, x') \in {\mathbb {S}}^d$, qui se concentrent sur un {\'e}quateur, montre que $p=2$ est le plus grand indice pour lequel cette propri{\'e}t{\'e} est vraie si on {\'e}tudie la norme $L^p$ de $\chi_\lambda f \,\chi_\mu g$. Le m{\^e}me calcul montre que~\eqref{eqbilin1f} est optimal sur  ${\mathbb S}^2$ tandis que sur $\mathbb{S}^d$ ($d\geq 3$), l'optimalit{\'e} (modulo la perte logarithmique) est obtenue en consid{\'e}rant les harmoniques sph{\'e}riques {\em zonales} qui se concentrent sur deux p{\^o}les de la sph{\`e}re.\par
D'autres applications du th{\'e}or{\`e}me~\ref{th2}  seront d{\'e}velopp{\'e}es dans un article ult{\'e}rieur.
\selectlanguage{english}
\section{Introduction}
\label{}



\selectlanguage{english}
Let $(M,g)$ be a compact smooth Riemannian manifold without boundary of dimension $d$  and 
$\Delta$
be the Laplace operator on functions on $M$. In~\cite{BGT7}, we proved a bilinear estimate for the spectral projectors of $\Delta$ in the case $d=2$. Our goal here is to extend this result to higher dimensions. 
\vskip 0.2 cm
\begin{theorem}\label{th1}
Let $\chi\in {\cal S}(\R)$. For $\lambda \in {\mathbb R}$,  denote by  
$\chi_\lambda = \chi(\sqrt{-\Delta}-\lambda)$ the spectral projector around $\lambda$.
 There exists $C$ such that for any $\lambda, \mu \geq 1 $, $f,g \in L^{2}(M)$,
\begin{equation}\label{eqbilin}
\|\chi_\lambda f \,\chi_\mu g
\|_{L^{2}(M)}
\leq C \Lambda(d, \min (\lambda, \mu)) \|f\|_{L^{2}(M)}\|g\|_{L^{2}(M)}
,
\end{equation}
{ with }
\begin{equation*}
\Lambda(d,\nu) = \begin{cases}
  \nu^{\frac{1}{4}} &\text{ if $d=2$}\\
  \nu^{\frac{d-2}{2}}\log^{1/2} (\nu) &\text{ if $d=3$}\\
  \nu^{\frac{d-2}{2}} &\text{ if $d\geq 4$} 
\end{cases}
\end{equation*}
Moreover for any  $\lambda, \mu, \nu \geq 1$, $f,g, h \in L^{2}(M)$, the following trilinear estimate holds
\begin{equation}\label{eqtrilin}
\|\chi_{\lambda}f\,  \chi_{\mu}g \, \chi_{\nu}h
\|_{L^{2}(M)}
\leq
C\left(\frac{ \lambda \mu \nu} {\max ( \lambda, \mu, \nu)}\right)^{\frac{2d-3}4}\|f\|_{L^{2}(M)}\|g\|_{L^{2}(M)}\|h\|_{L^{2}(M)}.
\end{equation}
\end{theorem}
Applying this result with  $f$ and $g$ two spherical harmonics on the sphere $\mathbb{S}^d$, we obtain
\vskip 0.2 cm
\begin{corollary}
\label{coeigen}
There exists $C>0$ such that if $H_p$ and $H_q$ are two spherical harmonics of respective degrees $p$ and $q$ greater than $1$,  
\begin{equation}\label{eqbilin11}
\|H_p H_q
\|_{L^{2}(\mathbb{S}^d)}
\leq C \Lambda(d, \min (p,q)) \|H_p\|_{L^{2}(\mathbb{S}^d)}\|H_q\|_{L^{2}(\mathbb{S}^d)}
\end{equation}
\end{corollary}
The {\em linear} versions of our theorem ($\lambda= \mu=\nu$), without the logarithmic loss for $d=3$, are due to Sogge~\cite{SO86,SO88,S}. In the case $d=2$ our proof in~\cite{BGT7} was inspired by H{\"o}rmander's work~\cite{Ho73} on Carleson Sj{\"o}lin type operators. The proof we present here is different even for $d=2$ and relies on  a bilinearization of the arguments in~\cite{SO86,SO88,S}. More precisely, after several reductions we reduce the matter to two (micro)-local {\em linear} estimates of quite a different nature. The fact that in~\eqref{eqbilin} the highest frequency disappears completely in the right hand side of the estimate is crucial in the applications to the non-linear Schr{\"o}dinger equation (see~\cite{BGT7}). A simple computation on the eigenfunctions $e_n(x_1, x_2,x') = (x_1+ix_2)^n$, $(x_1, x_2, x') \in {\mathbb {S}}^d$, which concentrate on the equator, shows that $p=2$ is the highest index for which this phenomenon occurs if one studies the $L^p$ norm of $\chi_\lambda f \, \chi_\mu g$. The same computation shows that~\eqref{eqbilin11} is optimal on ${\mathbb S}^2$. For $d\geq 3$, the optimality of~\eqref{eqbilin11} (except for the log loss) can be deduced by  considering the {\em zonal} eigenfunctions which concentrate on two poles of the sphere.\par
We close this introduction by mentioning that further applications of Theorem~\ref{th1} to the non linear Schr{\"o}dinger equation will be pursued in a forthcoming paper.
\section{Bilinear Estimates}
In this part we are going to give an outline of proof of the estimate~\eqref{eqbilin}. The proof  of~\eqref{eqtrilin} is similar. We assume $1\leq \lambda \leq \mu$. The estimate~\eqref{eqbilin} for any non trivial choice of $\chi$ implies~\eqref{eqbilin} for all $\chi$. Using a parametrix for the solution of the wave equations, Sogge~\cite{S} shows that for a suitable choice of the function $\chi$, we have:
$$
\chi_\lambda f
=
\lambda^{\frac {d-1} 2} T_{\lambda}f+R_{\lambda}f,
\qquad
\|R_{\lambda}f\|_{L^{\infty}}\leq C\|f\|_{L^2}
$$
 and in a coordinate system close to  $x_0\in M$,
$$
T_{\lambda}f(x)=\int_{\R^d}e^{i\lambda\varphi(x,y)}a(x,y, \lambda)f(y)dy
$$
where $a(x,y,\lambda)$ is a polynomial in $\lambda^{-1}$ with smooth coefficients supported in the set $$\{(x,y)\in {\mathbb R}^{2d}; \, |x| \leq \delta \ll \frac{\varepsilon}{C}\leq |y|\leq C\varepsilon\}
$$
 and   $-\varphi(x,y)=d_{g}(x,y)$ is the geodesic distance between $x$ and $y$.\par
 In  geodesic coordinates  $y=\exp_{0}(r\omega)$, $r >0$, $\omega \in \mathbb{S}^{d-1}$ we have
$$
T_{\lambda}f(x)=
\int_{0}^{\infty}\int_{\omega\in \mathbb{S}^{d-1}}
e^{i\lambda\varphi_{r}(x,\omega)}
a_r(x,\omega,\lambda)f_r(\omega)drd\omega \stackrel{\text{def}}{=} \int_{0}^{\infty}T_{\lambda}^{r}f_r (x) dr.
$$
where
$$ dy = \kappa(r,\omega) dr d\omega, \quad\varphi_r(x, \omega) = \varphi(x,r,\omega), \quad a_r(x, \omega,\lambda) = \kappa(r,\omega)a(x,r,\omega,\lambda), \quad f_r(\omega) = f(r,\omega).
$$
We get
$$
(T_{\lambda}f\,{T}_{\mu}g)(x)
=
\int_{\varepsilon/C}^{C\varepsilon}
\int_{\varepsilon/C}^{C\varepsilon}
T_{\lambda}^{r}f_r(x)\,
{T}_{\mu}^{q}g_q(x)drdq,
$$
and the Minkowski inequality shows that to prove Theorem~\ref{th1} it is enough to show, uniformly for $1\leq\lambda\leq \mu$,
\begin{equation}\label{1}
\|{T}_{\lambda}^{r}f\,{T}_{\mu}^{q}g\|_{L^2}
\leq
C\Lambda(d,\lambda)(\lambda \mu)^{-\frac{(d-1)}2}
\|f\|_{L^2}\|g\|_{L^2}.
\end{equation}
We have 
\begin{equation*}
{T}_{\lambda}^{r}f\,{T}_{\mu}^{q}g = \int_{\mathbb{S}^{d-1}\times \mathbb{S}^{d-1}}
e^{i\lambda\varphi_{r}(x,\omega)+i\mu\varphi_{q}(x,\omega')}
a_r(x,\omega,\lambda){a}_q(x,\omega',\mu){f}(\omega){g}(\omega')d\omega d\omega'.
\end{equation*}
 We first reduce the study to the case where $y\sim y'$ or $y\sim -y'$.
\vskip 0.2 cm
\begin{lemma}
\label{le1.2} For any $\omega'_0\in \mathbb{S}^{d-1}, \alpha>0$ there exist $c>0$ such that if $ |\omega_{0}- \omega_{0}'|\geq \alpha$ and $|\omega_{0}+\omega_{0}'|\geq \alpha$ then there exist $\varepsilon>0$ and a coordinate system in a neighbourhood of $\omega_{0}$, $\omega= (\omega_{1}, \dots , \omega_{d-1})$ such that if $\omega'$ is close to $\omega'_0$, $\omega$ close to $\omega_0$ and $\varepsilon/C \leq r,q \leq C\varepsilon$, then 
\begin{equation}\label{eq.nondeg}
 \text{det} \left[ \nabla _{x} \partial_{\omega_{1}} \varphi_{r}(x, \omega), \nabla_{x}\nabla_{\omega'}\varphi_{q}(x, \omega')\right] | \geq c >0.
\end{equation}
\end{lemma}
{\bf Proof.} Proceeding as in~\cite{BGT7} and using Gauss' lemma (see for example~\cite[3.70]{GaHuLa90}) we get 
\begin{equation}\label{eq1.4}
 \nabla_{x}\varphi_r(0,\omega)=\omega.
\end{equation}
Consequently the $(d-1)\times d$ matrix 
$ \nabla_{\omega}\nabla_{x}\varphi_r(0,\omega)
$
has rank $d-1$ and its range is the tangent plane to $\mathbb{ S}^{d-1}$ at $\omega$. Then, since we have supposed that $\omega_0\neq \pm \omega_0'$, the two tangent planes (seen as vector spaces) corresponding to $y$ and $y'$ are different and this implies the existence of the choice of the coordinate $\omega_{1}$ such that
$$\text{det} \left[ \nabla _{x} \partial_{\omega_{1}} \varphi_{r}(0, \omega_0), \nabla_{x}\nabla_{\omega'}\varphi_{r}(0, \omega'_0)\right] \neq 0$$
By continuity we get~\eqref{eq.nondeg} for $(x,r,\omega,\omega')$ close to $(0,q,\omega_0, \omega'_0)$. \qed \par
Consider now the operator (with frozen $\theta = (\omega_{2}, \dots , \omega_{d-1})$)
$$
{ {T}}_{\lambda}^{r,\theta}f(x)\,{{T}}_{\mu}^{q}g(x)
=
\int_{ {\mathbb R}\times {\mathbb R}^{d-1}} 
e^{i(\lambda\varphi_{r}(x,\omega_1, \theta)+\mu\varphi_{q}(x,\omega'))}
a_{r}(x,\omega_1, \theta,\lambda){a}_{q}(x,\omega',\mu)f(\omega_1)g(\omega')d\omega_1 d\omega'.
$$
We can write
$$
\|{ {T}}_{\lambda}^{r,\theta}f\,{{T}}_{\mu}^{q}g\|_{L^2}^2
=
\int
K(\omega_1,\sigma_1,\omega',\sigma')
f(\omega_1)g(\omega')\overline{f(\sigma_1)g(\sigma)}
d\omega_1 d\sigma_1 d\omega' d\sigma'
$$
where, according to Lemma~\ref{le1.2}, by integrations by parts in $x$, we get easily
$$
|K(\omega_1,\sigma_1,\omega',\sigma')|
\leq
C_{N}(1+ \lambda|\omega_1- \sigma_1|+\mu|\omega'- \sigma'|)^{-N}.
$$
Schur's Lemma gives
\begin{equation*}
\begin{aligned}
\|{ {T}}_{\lambda}^{r,\theta}f\,{{T}}_{\mu}^{q}g\|_{L^2}^2
&\leq C_N \left[ \int \frac{dtd\varrho}{(1+\lambda|t|+\mu|\varrho|)^N}
\right]\|f\|_{L^2}^{2}\|g\|_{L^2}^{2}
&\leq 
C{\mu^{-{(d-1)} }}\lambda^{-1}\|f\|^2_{L^2}\|g\|^2_{L^2}.
\end{aligned}
\end{equation*}
Using Minkowski inequality, this implies
\vskip 0.2 cm
\begin{proposition}\label{prop.3}
Suppose that $a_r(x,\omega,\lambda) a_q(x,\omega', \mu)$ is supported in a set where 
$$ |\omega- \omega'|\geq \alpha \text{ and } |\omega+\omega'|\geq \alpha$$ and that $q$ and $r$ are close. Then
$$
\|{ {T}}_{\lambda}^{r}f\,{{T}}_{\mu}^{q}g\|_{L^2}
\leq 
{C}{\mu^{-\frac{d-1} 2}}\lambda^{-\frac 1 2}\|f\|_{L^2}\|g\|_{L^2}.
$$
\end{proposition}
Remark that Proposition~\ref{prop.3} gives (except for the log loss in dimension $3$) the estimate~\eqref{1} if $d\geq 3$ and is better if $d=2$.\par
We are now left with the following two cases:
\begin{enumerate}
\item\label{casei}
$a_r$, ${a}_q$ are localized close to $\omega= \omega_{0}$,
$\omega'=\omega_0$ respectively.
\item\label{caseii}
$a_r$, ${a}_q$ are localized close to $\omega=\omega_0$,
$\omega'=-\omega_0$ respectively.
\end{enumerate}
We are going to study the case \eqref{casei}, the case \eqref{caseii} being  similar. We follow Sogge's strategy~\cite{S}, in the spirit of the arguments of Stein~\cite{St93}.
\vskip 0.2 cm
\begin{lemma}
\label{le1.3}
The phase $\varphi_{q}(x,\omega)$
is a Carleson Sj\"olin phase: near any point $(x_{0}, \omega_{0})$, one can choose a splitting of the variable $x= (t,z)$
such that
\begin{enumerate}
\item \label{case1}For fixed $t$ the phase $\varphi_{q}(t,z,\omega)$ is uniformly non degenerate:
\begin{equation}\label{eq1.0} \left|\text{det}\left(\frac{\partial^2 \varphi_{q}(t,z,\omega)} { \partial z_{j}\partial{\omega_{i}} }\right)\right|\geq c >0
\end{equation}
\item \label{case2}Let 
$$S_{t,z}= \{ \nabla_{t,z} \varphi_q (t,z,\omega), \, \omega \sim \omega_0\}$$
Then $S$ is according to~\eqref{case1} a smooth hypersurface in ${\mathbb R}^d$ and  it has non-vanishing principal curvatures: 
denote by $\nu (t,z,\omega)$ the normal unit vector to the surface $S$ at the point $\varphi_q (t,z,\omega)$. Then
\begin{equation}\label{eq1.2}
\left|\text{det}\langle\frac{\partial^2 } { \partial \omega_{j}\partial{\omega_{i}} } \nabla_{t,z}\varphi_{q}(t,z,\omega), \nu( t,z,\omega)\rangle\right|\geq c> 0 
\end{equation}
\end{enumerate}
Furthermore if $r$ is close to $q$ and $\omega_{0}$ close to $\omega'_{0}$ or close to $-\omega'_{0}$, then we can choose the same splitting for the phases $\varphi_{q}(x,\omega)$ and $\varphi_{r}(x, \omega')$. 
\end{lemma}\vskip 0.2 cm
{\bf Proof.} We may assume $\omega_0=(1,0,\dots, 0)$. We choose $t=x_1, z= (x_2, \dots, x_d)$. Equations~\eqref{eq1.0} and~\eqref{eq1.2} at the point $t=0, z=0$ are an easy consequence of~\eqref{eq1.4} and the lemma follows by continuity.\qed \par \vskip 0.2 cm 
Using~\eqref{eq1.0} we deduce the following refinement of the $L^2$ boundedness of the spectral projector.
\vskip 0.2 cm 
\begin{proposition}
\label{prop.2}
The operator
$$g\in L^2(M) \mapsto T_{\nu}^{r}g(t,z)\in L^\infty( {\mathbb R}_{t}; L^2( {\mathbb R}^{d-1}_{z}))$$
is continuous with norm bounded by $C \nu^{ -(d-1)/2}$.
\end{proposition}
\vskip 0.2 cm 
 On the other hand, using the dispersion property~\eqref{case2} in Lemma~\ref{le1.3} leads to the following. 
\vskip 0.2 cm
\begin{proposition}
\label{prop.1}In the coordinate system of Lemma~\ref{le1.3}, we have:
\begin{equation*}
\| T_{\nu}^{r}f(t,z)\|_{L^2( {\mathbb R}_{t}; L^\infty( {\mathbb R}^{d-1}_{z}))}\leq C\Lambda(d, \nu) \nu^{-(d-1)/2} \| f\|_{L^2}.
\end{equation*}
\end{proposition}
Before giving the proof of this result, let us show how to finish the proof of Theorem~\ref{th1}. Putting together  Propositions~\ref{prop.1} and~\ref{prop.2} we get (using that the same splitting can be chosen for $\varphi_{r}$ and $\varphi_{q}$) that 
$$(f,g)\in L^2_\omega \times L^ 2_{\omega'} \mapsto T_{\lambda}^{r}f(t,z)\times T_{\mu}^{q}g(t,z)\in L^2( {\mathbb R}_{t}; L^2( {\mathbb R}^{d-1}_{z}))$$
is continuous with norm bounded by 
$$
C (\lambda \mu)^{-(d-1)/2} \Lambda(d, \lambda)
$$
which is~\eqref{1}\qed \par

Finally let us come to the proof of Proposition~\ref{prop.1}. By a $TT^*$ argument, it is enough to estimate the norm of
$$T_{\nu}^{r}(T_{\nu}^{r})^*: L^{2}( {\mathbb R}_{t'}; L^1( 
{\mathbb R}^{d-1}_{z'}))\mapsto L^2( {\mathbb R}_{t}; L^\infty( 
{\mathbb R}^{d-1}_{z})).$$
But the kernel of this operator is
$$K(t,z,t',z')= \int e^{i\nu( \varphi_{r}(t,z,\omega)- \varphi_{r}(t',z', \omega))}a_r(t,z,\omega,\nu) a_r(t',z', \omega, \nu)d \omega$$
and using the second part in Lemma~\ref{le1.3} we can show that
$$| K(t,z,t',z')|\leq \frac { C} { \left(1+ \nu |(t,z)- (t',z')|\right)^{(d-1)/2}}\leq \frac { C} {\left( 1+ \nu |t-t'|\right)^{(d-1)/2}} $$
Indeed,
$$\nabla_\omega \langle \nabla _{t,z}\varphi_r(t,z,\omega), n\rangle =0 \Leftrightarrow n= \nu(t,z,\omega)$$
and  taking into account that
\begin{equation*}
\varphi_{r}(t,z,\omega)- \varphi_{r}(t',z', \omega)
= \langle\nabla_{t,z}\varphi_{r}(t,z,\omega), (t,z)- (t',z')\rangle + {\mathcal {O}}( \| (t,z)- (t',z')\|^2),\end{equation*}
 we see that if $(t,z)- (t', z')$ is in a small conic neighbourhood of the critical direction $\nu(t,z)$ then, in view of~\eqref{eq1.2}, we can apply the stationnary phase formula to the integral in $\omega$. Otherwise, we can integrate by parts in $\omega$ and using that
$$|\nabla_{\omega}\left(\varphi_{r}(t,z,\omega)- \varphi_{r}(t',z', \omega)\right)| \geq c |(t,z)- (t',z')|
$$
we get for any $N\in {\mathbb N}$,
$$| K(t,z,t',z')|\leq \frac { C_N} { (1+ \nu |(t,z)- (t',z')|)^N}
$$
 which is better (see~\cite[p 63]{S}).\par
 We conclude using  the classical one dimensional Young inequality
\begin{equation*}
\|T_{\nu}^{r}(T_{\nu}^{r})^*\|_{\mathcal {L}\left(L^2(]-1,1[_{t'}; L^1_{z'});L^2(]-1,1[_t; L^\infty_{z})\right)} \leq C \int_{|s|\leq 2} \frac {ds} {(1+\nu |s|)^{ (d-1)/2}}
\leq \begin{cases} C \nu^{-1/2} & \text{ if } d=2\\
C \nu^{-1}\log (\nu) & \text{ if } d=3\\
C \nu^{-1} & \text{ if } d\geq 4.
\end{cases}
\end{equation*}

\end{document}